\documentclass[11pt]{article}

\usepackage{amsmath}
\usepackage{amssymb}

\numberwithin{equation}{section}

\newtheorem{thm}{Theorem}[section]

\newtheorem{prop}[thm]{Proposition}

\newtheorem{rem}[thm]{Remark}

\newcommand{\grad}{\operatorname{grad}}
\newcommand{\supp}{\operatorname{supp}}
\newcommand{\stab}{\operatorname{stab}}

\def\mod{\text{Mod}_{g,n}}
\def\ga{\Gamma}

\def\gag{\Gamma_{g,n}}

\def\R{\mathbb R}

\def\Z{\mathbb Z}

\def\T{\mathcal T}

\def\E{\underline{E}\Gamma_{g,n}}

\title{A cofinite universal space for proper actions for mapping class groups\footnote{2000 Mathematics Subject Classification Primary: 57M07; Secondary: 30F60, 20F34.}}
\author{Lizhen Ji
\\ Department of Mathematics\\ University of Michigan\\ Ann Arbor, MI 48109
\\ \&  \\ Scott A. Wolpert\\ Department of Mathematics\\ University of Maryland
\\ College Park, MD 20742}

\date{January 5, 2009}
\begin{document}

\maketitle

\begin{abstract}
We prove that the mapping class group $\gag$ for surfaces of negative
Euler characteristic has a cofinite universal space $\E$ for proper actions (the resulting quotient is a finite $CW$-complex).  The approach is to construct a truncated Teichm\"{u}ller space $\T_{g,n}(\varepsilon)$  by introducing a lower bound for the length of shortest closed geodesics and showing that  $\T_{g,n}(\varepsilon)$ is a $\gag$ equivariant deformation retract of the Teichm\"{u}ller space $\T_{g, n}$.  The existence of such a cofinite universal space is important in the study of the cohomology of the group $\gag$. As an application, we note that there are only finitely many conjugacy classes of finite subgroups of $\gag$. 
Another application is that the rational Novikov conjecture in K-theory
holds for $\gag$.
\end{abstract}

\section{Introduction}
The mapping class group $\gag$ for surfaces of genus $g$ with $n$ punctures 
is a natural group and its structure has been intensively
studied by many authors employing a range of methods. A natural question
concerns its classifying space.
Since $\gag$ contains torsion-elements, the group
cannot admit a finite dimensional classifying space $B\mod$.
Ivanov [Iv1] [Iv2] did prove that for any torsion-free subgroup $\ga$ of $\gag$ of finite
index there exists a $B\ga$-space given by a finite $CW$-complex.  Given the foundational Deligne-Mumford  quasiprojectivity of $\gag\backslash \T_{g,n}$, a proof for torsion-free subgroups also follows from the original result of Lojasiewicz [Lo].

Recall that for a discrete group $\ga$, the universal covering space
$E\ga=\widetilde{B\ga}$ of $B\ga$ is a universal space for proper and fixed-point free actions of $\ga$.   For groups $\ga$ containing torsion elements, a closely related space
 is the universal space $\underline{E}\ga$ for proper actions of $\ga$, which is unique up to homotopy.
 Briefly,  for any discrete group $\ga$, a $CW$-complex $E$
is {\em a universal space for proper actions} of $\ga$ if the following conditions are satisfied:
\begin{enumerate}
\item $\ga$ acts properly on $E$ and hence the stabilizer in $\ga$ of every point is finite.
\item For any finite subgroup $H\subset \ga$, the set $E^{H}$ of fixed points of $H$
is nonempty and contractible. In particular, $E$ is contractible.
\end{enumerate}
 If the quotient $\ga\backslash E$ is a finite $CW$-complex, the quotient is called {\em a cofinite universal space  for proper actions}. The existence of such a cofinite $\ga$-space is important in studying
the cohomology of $\ga$ and the Novikov, Baum-Connes and related conjectures for $\ga$
(see [BCH] and [Lu] for systematic discussions).
Applying the results of [Brow1] and [BHM], it can be seen easily that the existence
of a cofinite universal space $\underline{E}\ga$ implies that the rational Novikov conjecture
in K-theory holds for $\ga$ (see Proposition \ref{r-Novikov} below).

It is known and can be shown easily that the Teichm\"{u}ller space $\T_{g, n}$ is a finite dimensional, non-cofinite $\E$-space (see [Lu] and Proposition \ref{finite-dim} below). 
As mentioned above, it is important for various purposes to find a cofinite universal space
(a space with quotient a finite $CW$-complex) and the question on the existence of a cofinite $\E$ has been raised by several authors (see [Lu]). 
The main result of this note is a positive answer (see Theorems \ref{main} and \ref{contractible} below).  When $n\geq 1$, the result also follows from an explicit triangulation of $\T_{g,n}$ and an associated retraction onto a cofinite subspace, see Harer 
[Har]. After the results of this paper were obtained,
we found that the case of $n=0$ was also established already in [Brou1] and [Mi].
On the other hand, the method here is different and works in both cases.

Specifically, let $S_{g,n}$ be a surface of genus $g$ with $n$ punctures. Assume that $2g+n>2$.
Let $\T_{g, n}$ be the associated Teichm\"{u}ller space of homotopy marked hyperbolic metrics on $S_{g,n}$.  It is well-known that $\T_{g,n}$ is diffeomorphic to $\R^{6g-6+2n}$ and hence
contractible, and $\gag$ acts properly on $\T_{g, n}$. 
Using the solution of the Nielsen realization problem and the Weil-Petersson metric for $\T_{g,n}$, it can be shown that $\T_{g, n}$ is the universal space for proper actions of $\T_{g,n}$
(see [Lu, \S 4.10] for discussion of this result and Proposition \ref{finite-dim} below for a detailed proof).

It is also well-known that the quotient $\gag\backslash \T_{g,n}$ is noncompact,
hence $\T_{g,n}$ is not a cofinite $CW$-space.
In fact, the quotient is the moduli space of Riemann surfaces of type $(g,n)$.
Since Riemann surfaces can degenerate, their moduli space is noncompact. 
There are two basic ways to address the issue of noncompactness.

One way is to compactify $\gag \backslash \T_{g,n}$, or rather to partially
compactify  $\T_{g, n}$ to a manifold with corners 
$\overline{\T}^{BS}_{g,n}$ such that $\gag$ acts
properly with a compact quotient $\gag \backslash \overline{\T}^{BS}_{g,n}$.
This approach is analogous to considering the Borel-Serre compactification of locally symmetric spaces following [BS]. Such a partial compactification of $\T_{g,n}$ as a real analytic manifold
with corners was announced by Harvey [Hv1]. A proof and construction of
$\overline{\T}^{BS}_{g,n}$ as a smooth manifold with corners was later
given by Ivanov [Iv3].  The next step is to show that 
 the set $(\overline{\T}^{BS}_{g,n})^{H}$ of fixed points of a finite subgroup $H\subset \gag$
 is contractible. 
 For arithmetic subgroups $\ga$ of semisimple Lie groups,  a construction
 of $\underline{E}\ga$ in terms of the partial Borel-Serre compactification of the associated
 symmetric spaces was carried out in [Ji], where the geodesic action of parabolic
 subgroups and its connection to the boundary components was used crucially.
 Since it might be involved to carry out the second step for $\ga_{g,n}$, we 
introduce a different model for the classifying space.

Let $S\in \T_{g,n}$ be a marked surface with hyperbolic metric.
For every  closed curve $\sigma\subset S$, let $\ell_{S}(\sigma)$ be the length
of the unique geodesic in the free homotopy class of $\sigma$ with respect to the hyperbolic metric of $S$.
For a sufficiently small $\varepsilon >0$, define
\begin{equation}
\T_{g,n}(\varepsilon)=\{S \in \T_{g,n}\mid \text{ \ for every closed curve\ }
\sigma \subset S,\  \ell_{S}(\sigma)\geq \varepsilon\}.
\end{equation}

This is called the {\em truncated Teichm\"{u}ller space}.  The following is well-known.

\begin{prop}\label{cw-complex}
The truncated Teichm\"{u}ller space $\T_{g, n}(\varepsilon)$ is a real analytic manifold with corners and as a subset of $\T_{g,n}$ is invariant under the action of $\gag$.   The quotient $\gag\backslash \T_{g, n}(\varepsilon)$ is a compact real analytic orbifold with corners. In particular, $\T_{g, n}(\varepsilon)$  has the structure of a cofinite $\gag$-$CW$-complex.
\end{prop}

Note that the last statement on the structure of a cofinite $\gag$-$CW$-complex follows
from the existence of an equivariant triangulation in [Il]. 
For any torsion-free subgroup $\ga\subset \gag$, the quotient $\ga\backslash \T_{g,n}$ is a real analytic manifold with corners. The proof of Proposition \ref{cw-complex} will be given in \S 2 below.

The main result of the paper is the following.

\begin{thm}\label{main}
With the above notation, when $\varepsilon$ is sufficiently small,  there is an $\gag$-equivariant deformation retraction from $\T_{g, n}$ to $\T_{g, n}(\varepsilon)$.
\end{thm}

The corresponding moduli space retraction from $\gag\backslash\T_{g, n}$ to $\gag\backslash\T_{g, n}(\varepsilon)$ is suggestive of the realizations of the Borel-Serre compactifications for locally symmetric spaces [Gr] [Le1-2] [Sa].   
This result was stated in [Har, \S 3, iv)] but the proof is incomplete
(see \S 3 below for discussion).  One present purpose is to complete the considerations of [Har]. 
 It is also noteworthy to point out that for $\ga\subset\ga_{g,n}$ a torsion-free
 finite index subgroup, it was proved in [Iv1] [Iv2] that $\T_{g,n}(\varepsilon)$ is 
 a $\ga$-equivariant deformation retract. The assumption that $\ga$ is torsion-free
 was used in a essential way.

An immediate corollary of the theorem is the following.

\begin{thm}\label{contractible}
The truncated space  $\T_{g, n}(\varepsilon)$ is contractible.
For a finite subgroup $H\subset \gag$, the set of fixed points
$(\T_{g, n}(\varepsilon))^{H}$ is also contractible, and hence $\T_{g, n}(\varepsilon)$ is a cofinite universal space for proper actions of $\gag$.
\end{thm}

\begin{rem}
{\em After this result was proved, we learned of and received a copy
of the preprint [Mi] from L\"{u}ck.  The preprint discusses the Broughton considerations of [Brou1] for the existence of a cofinite classifying space of $\ga_g$  by introducing a 
suitable subspace of $\T_g$.  
It might be helpful to point out that the present proof is closely related to the corresponding result for arithmetic groups and works for all choices of $(g, n)$, while the method
in [Brou1] only works for the case $(g,0)$, since the Satake compactification
for the moduli space $\ga_g\backslash \T_g$ in [Ba] is only defined for the case $(g,0)$
and is used in an essential way.
}
\end{rem}

An interesting corollary of the above theorem is the following.
\begin{prop}\label{finite-many}
There are only finitely many conjugacy classes of finite subgroups of $\gag$.
\end{prop}
\begin{proof}
It is a general and known fact that if a group $\ga$ admits a cofinite universal space $\underline{E}\ga$
for proper actions, there are only finitely many conjugacy classes of finite subgroups.
It can be proved as follows.
Let $K$ be a compact subset of $\underline{E}\ga$ which is mapped surjectively
onto the quotient $\ga\backslash \underline{E}\ga$.
Then any finite subgroup of $\ga$ is conjugate to a finite subgroup
which has a fixed point in $K$.
Since there are only finitely many elements $\gamma\in \ga$ such that $\gamma K\cap K\neq\emptyset$,
it follows that there are only finitely many conjugacy classes of finite subgroups. 
\end{proof}

Since this result is of independent interest, we provide a second proof
in Proposition \ref{finite-many-1} below.  We note that the finiteness result is similar to the fact that any arithmetic subgroup of a Lie group has only finitely many conjugacy classes of finite subgroups (see [Se]).

\begin{rem}
{\em After the result in Proposition \ref{finite-many}  was proved, we learned that 
it was proved already by Bridson in [Bri, Theorem 6] using the above result of Harer  [Har, \S 3].
Since the proof of [Har, \S 3] is incomplete, we now complete the proof of [Bri].  An independent proof of Proposition \ref{finite-many} is given in the book manuscript of Farb and Margalit [FaM].
}
\end{rem}

\begin{rem}
{\em
It is interesting to note that [Brou2] first proved Proposition \ref{finite-many}
and then used it to construct a cofinite universal space for $\ga_{g,n}$.
Proposition \ref{finite-many} was proved in [Brou2] by showing there are 
only finitely many equivalence classes of finite orientation-preserving
topological actions on surfaces of a fixed genus. 
}
\end{rem}

Another corollary of Theorem \ref{contractible} is the following.

\begin{prop}\label{r-Novikov}
The group $\gag$ is of type $FP_\infty$ and hence in every degree
$H_*(\gag, \Z)$ is finitely generated, which implies that the rational
Novikov conjecture in K-theory holds for $\gag$.
\end{prop}
\begin{proof}
Since $\gag$ admits a cofinite universal space for proper actions and finite
groups are of type $FP_\infty$, the first statement on finite generation of $H_*(\gag, \Z)$ follows from [Brow1, Proposition 1.1], and the second statement follows from
the main result of [BHM] which states that the rational Novikov conjecture in K-theory holds
for any groups $\ga$ such that $H_*(\ga, \Z)$ is finitely generated in every degree.
\end{proof}

\begin{rem}
{\em 
We note that the original Novikov conjecture, i.e., the rational Novikov conjecture
in L-theory for $\gag$ was proved by Kida [Ki] and Hamenst\"{a}dt [Ham].
}
 \end{rem}
 
 \begin{rem}
 {\em
 After Proposition \ref{r-Novikov} was proved as above, we realized that it can also be proved
 as follows. It is known that $\gag$ admits a torsion-free subgroup $\ga'$ of finite index (see [Iv1] for a demonstration of this not-totally obvious result. See also [Har]).
 By [Iv1], $\T_{g,n}$-admits a $\ga'$-equivariant deformation retraction to
 a $\ga'$-compact subspace (a manifold with corners).  It follows that $\ga'$ admits
 a finite classifying space $B\ga'$, which in turn implies that $\ga'$ is of type $FP_\infty$.
 Then by [Brow2, Proposition 5.1, p. 197], $\gag$ is also of type $FP_\infty$,
and  the main result of [BHM] implies the second statement as above.  After Proposition \ref{r-Novikov} was proved, we learned of Storm's [St] note for combining results of Hamenst\"{a}dt and Kato on combable groups to also provide a proof.  
 }
 \end{rem}
 
\vspace{.1in}

\section{Mapping class groups and Teichm\"{u}ller spaces} 

First we recall definitions.  Let $S_{g,n}$ be an orientable surface of genus $g$ with $n$ punctures.  Let $\text{Diff}(S_{g,n})$ be 
the group of all diffeomorphisms of $S_{g,n}$ and $\text{Diff}^{+}(S_{g,n})$ the subgroup of orientation preserving diffeomorphisms.
The identity component $\text{Diff}^0(S_{g,n})$ is a normal subgroup of both
$\text{Diff}(S_{g,n})$ and $\text{Diff}^{+}(S_{g,n})$.  The quotient 
$\pi_{0}(\text{Diff}^{+}(S_{g,n}))=\text{Diff}^{+}(S_{g,n})/\text{Diff}^0(S_{g,n})$ 
is called the {\em mapping class group} of type
$(g,n)$, denoted by $\ga_{g,n}$, and $\pi_{0}(\text{Diff}(S_{g,n}))
=\text{Diff}(S_{g,n})/\text{Diff}^0(S_{g,n})$
is called the {\em extended mapping class group}, denoted by $\tilde{\ga}_{g,n}$.
Clearly, $\ga_{g,n}$ is an index two subgroup of $\tilde{\ga}_{g,n}$.
 
 The mapping class group $\ga_{g,n}$ is closely related to arithmetic groups.
 In fact, for $S_{g,n}$ the torus $\R^{2}/\Z^{2}$, i.e., $g=1$ and $n=0$,
 $\ga_{1}$ is canonically isomorphic to $SL(2, \Z)$. 

An important step for understanding $\ga_{g,n}$ is to consider its action
on the Teichm\"{u}ller space $\T_{g,n}$ associated with $S_{g,n}$.
Assume that the Euler number of $S_{g,n}$ is negative, i.e., $2g-2+n>0$.
The surface $S_{g,n}$ admits complete hyperbolic metrics of finite volume.
Recall that a marked hyperbolic metric on $S_{g,n}$ is a surface $S$ with a complete
hyperbolic metric $ds$ together with a diffeomorphism $\varphi:S_{g,n}\to S$.
Two marked hyperbolic metrics $(S_{1}, ds_{1}, \varphi_{1})$
and $(S_{2}, ds_{2}, \varphi_{2})$ are {\em equivalent}
if there exists an isometry $\phi:S_{1}\to S_{2}$ such that
$\phi\circ \varphi_{1}:S_{g,n}\to S_{2}$ is isotopic to
$\varphi_{2}:S_{g,n}\to S_{2}$.  
The set of equivalence classes of marked hyperbolic metrics on $S_{g,n}$
forms the {\em Teichm\"{u}ller space} $\T_{g, n}$.

\begin{prop}  
The groups $\ga_{g, n}$ and $\widetilde{\ga}_{g,n}$ act properly discontinuously on $\T_{g,n}$.  
The quotient $\ga_{g, n}\backslash \T_{g,n}$ is the moduli space of
complex structures on $S_{g,n}$ of type $(g,n)$ and $\widetilde{\ga}_{g,n}\backslash \T_{g,n}$ is the moduli space of complete hyperbolic metrics on $S_{g,n}$ of type $(g,n)$.
\end{prop}

It is well-known that $\T_{g,n}$ is diffeomorphic to $\R^{6g-6+2n}$.
One way to display this is to introduce Fenchel-Nielsen coordinates.   
Let $\Sigma=\{\sigma_{1}, \cdots, \sigma_{d}\}$, $d=3g-3+n$, be a maximal collection
of disjoint simple closed curves on $S_{g,n}$, no curve freely homotopic to a point
or a puncture, and no pair of curves freely homotopic.  
The complement in $S_{g,n}$ of $\Sigma$ consists of pairs of pants, subsurfaces homeomorphic to the complement of three points in a sphere.  
The collection $\Sigma$ is called a {\em pants decomposition}.
Then for each point $(S, ds, \varphi)\in \T_{g,n}$, each curve $\varphi(\sigma_{i})$,
$i=1, \cdots, d$, contains a unique simple closed geodesic in its free homotopy class with respect
to the metric $ds$.
Denote the length of the geodesic by $\ell_{S}(\sigma_{i})$.
Once a base point $S_{1}\in \T_{g,n}$ is chosen, there are also associated twist parameters
$\theta_{S}(\sigma_{i})$.
Together, the data provides the Fenchel-Nielsen coordinates:
\begin{equation}
\pi_{\Sigma}:T_{g,n}\to \R^{6g-6+2n};\quad (S, ds,\varphi)\mapsto (\ell_{S}(\sigma_{1}), \theta_{S}(\sigma_{1}); \cdots;
\ell_{S}(\sigma_{d}), \theta_{S}(\sigma_{d})).
\end{equation}
It is known that $\T_{g, n}$ is a complex manifold, and the map $\pi_{\Sigma}$ is a real analytic diffeomorphism.

As a complex manifold, $\T_{g,n}$ admits a  K\"{a}hler metric, the Weil-Petersson metric.
At a point $S\in \T_{g, n}$, the dual of the holomorphic tangent space is canonically
identified with the vector space $Q(S)$ of holomorphic quadratic differentials.
The Hermitian product on $Q(S)$,
\begin{equation}
\langle \varphi_1, \varphi_2\rangle=\int_S \varphi_1 \overline{\varphi}_2 ds^{-2},
\end{equation}
where $ds^2$ is the area form of $S$, defines the Weil-Petersson metric on $\T_{g,n}$,
 denoted by $ds_{WP}$.

The following properties of the metric $ds_{WP}$ are known (see [Wo1-3]).

\begin{prop}\label{w-p-metric}
\begin{enumerate}
\item $\tilde{\ga}_{g,n}$ acts isometrically with respect to $ds_{WP}$.
\item The sectional curvatures of $ds_{WP}$ are negative.
\item $(\T_{g,n}, ds_{WP})$ is convex in the sense that each pair of points is connected
by a unique geodesic. 
\item For a finite subgroup $H$ of $\tilde{\ga}_{g,n}$ or $\ga_{g,n}$, a finite $H$-invariant compact subset of $\T_{g,n}$ has compact convex hull, and hence every finite subgroup of $\tilde{\ga}_{g,n}$ has a fixed point in $\T_{g,n}$.
\end{enumerate}
\end{prop}

Using these properties, we prove the following refinement of the Nielsen realization theorem.

\begin{prop}\label{finite-dim}
The space $\T_{g,n}$ is a universal space for proper actions of $\tilde{\ga}_{g,n}$ and $\gag$.
\end{prop}
\begin{proof} We consider only $\tilde{\ga}_{g,n}$.
Since $\T_{g,n}$ is a real analytic manifold and $\tilde{\ga}_{g,n}$ acts real analytically and properly on $\T_{g, n}$, the existence of an equivariant triangulation 
in [Il] implies that $\T_{g,n}$ is a proper $\tilde{\ga}_{g,n}$-$CW$-complex.
We need to show that for a finite subgroup $H\subset \tilde{\ga}_{g,n}$,
the fix-point set $\T_{g,n}^{H}$ is nonempty and contractible. 

By the solution of the Nielsen realization problem in [Ke] (see also [Wo2]), $\T_{g,n}^{H}$ is nonempty.  By Proposition \ref{w-p-metric}, $\T_{g,n}$ with the Weil-Petersson metric is
negatively curved and convex.   Since $\tilde{\ga}_{g,n}$ acts isometrically on $\T_{g,n}$, the set $\T_{g,n}^{H}$ is a totally geodesic submanifold, which is in turn also convex and hence contractible.
\end{proof}

\begin{rem}
{\em For the orientation preserving case, there is a different approach 
as indicated in the discussion of [Lu, \S 4.10].
Kerckhoff [Ke] considers the important result of Thurston that each pair of points $x, y$ of $\T_{g,n}$ is connected by a unique left earthquake, denoted by $\widetilde{x, y}$.
Note that the notion of the {\em left} earthquake depends on the orientation
of the underlying surface, and the ordering of $x$ and $y$ is also important. In fact, the earthquake from $x$ to $y$, $\widetilde{x, y}$, differs from the earthquake from $y$
to $x$, $\widetilde{y,x}$.
Earthquake paths are natural; for an (orientation preserving) element $h\in \gag$,
$h$ maps $\widetilde{x, y}$ to $\widetilde{h\cdot x, h\cdot y}$.
Naturality provides that if $H$ is a finite subgroup of $\gag$,
then for a pair $x, y\in \T_{g,n}^{H}$, $\widetilde{x, y}$
is contained in $\T_{g,n}^{H}$, and hence the fixed-point set $\T_{g,n}^{H}$ is contractible.
Naturality further provides that the left earthquake exponential map is
$\gag$-equivariant.

On the other hand,  if elements of $H\subset \tilde{\ga}_{g,n}$ are orientation reversing,
then $\widetilde{x, y}$ may not be contained in $\T_{g,n}^{H}$.
An example begins with hyperbolic surfaces with topologically conjugate mirror symmetry with respect to a separating simple closed geodesics.
}
\end{rem}

\begin{rem}
{\em In the above proof of Proposition \ref{finite-dim}, 
we used the unique Weil-Petersson metric geodesic connecting $x$
and $y$, denoted by $\overline{x, y}$. Since $\tilde{\ga}_{g,n}$ acts isometrically on 
$\T_{g,n}$, it follows for any subgroup $H\subset \tilde{\ga}_{g,n}$,
$\overline{x, y}\subset \T_{g,n}^{H}$ whenever
$x, y\in \T_{g,n}^{H}$, and hence the fixed-point set $\T_{g,n}^{H}$ is contractible.  
The proof can be made for any $\tilde{\ga}_{g,n}$ invariant metric with unique geodesics connecting pairs of points.  
}
\end{rem}

The quotient space $\ga_{g,n}\backslash \T_{g,n}$ is noncompact.
In particular, $\T_{g,n}$ is not a cofinite $\underline{E}\gag$-space; the quotient space is not a finite $CW$-complex.  
One way to understand the structure near infinity of $\ga_{g,n}\backslash \T_{g,n}$
is through the notion of {\em rough fundamental domains} and {\em Bers regions}.
For a pants decomposition $\Sigma$, constants $C>0$ and $\theta_{0}>0$,
define
\begin{equation}
\mathcal B_{\Sigma}=\mathcal B_\Sigma(C, \theta_{0})=\pi_{\Sigma}^{-1}
(\{(\ell_{1}, \theta_{1};\cdots; \ell_{d}, \theta_{d})\in \R^{6g-6+2n}\mid
0< \ell_{i}\leq C,\, |\theta_{i}|\leq \theta_0,\, i=1, \cdots, d\}).
\end{equation}
The region is called the {\em Bers region associated with $\Sigma$}.

The following result of Bers is analogous to the reduction theory for arithmetic groups.
See [Bu, Theorems 5.1.2  and 6.6.5] for proofs.

\begin{prop}\label{bers-region}
Up to isotopy, there are only finitely many pants decompositions $\Sigma_{1}, \cdots, \Sigma_{m}$
of $S_{g,n}$.  For each $\Sigma_{i}$, there is a Bers region $\mathcal B_{\Sigma_{i}}$ such that the projection map $\pi:\T_{g,n}\to \ga_{g,n}\backslash \T_{g,n}$
restricts to a finite-to-one map on $\mathcal B_{\Sigma_{i}}$.   The images of the Bers regions
$\pi(\mathcal B_{\Sigma_i})$ cover the quotient $\ga_{g,n}\backslash \T_{g,n}$.
\end{prop}

\noindent{\em Proof of Proposition \ref{cw-complex}}.  
We first observe that since  $\T_{g,n}(\varepsilon)$ is defined by a condition for all simple closed curves, the subset is invariant under $\gag$.  
We next observe that the $\gag$ quotient is compact.
In particular, for a pants decomposition $\Sigma=\{\sigma_{1}, \cdots, \sigma_{d}\}$ 
there is the associated Fenchel-Nielsen coordinate system:
\begin{equation}
\T_{g, n}\to \R^{6g-6+2n}, \quad (S,ds,\varphi) \mapsto (\ell_{S}(\sigma_{1}),
\theta_{S}(\sigma_{1});\cdots; \ell_{S}(\sigma_{d}),
\theta_{S}(\sigma_{d})).
\end{equation}
For $\varepsilon>0$, $\theta_{0}>0$ and $C>0$, introduce the {\em truncated Bers region}: 
\begin{equation}
\mathcal B_{\Sigma}(\varepsilon)=\{S\in T_{g,n}\mid
C\geq\ell_{S}(\sigma_{i})\geq \varepsilon,\, |\theta_{S}(\sigma_{1})|\leq \theta_0,\,  i=1, \cdots, d\}.
\end{equation}
The truncated regions are compact submanifolds with corners.
By Proposition \ref{bers-region}, the images of  
 finitely many truncated Bers regions 
$\mathcal B_{\Sigma}(\varepsilon)$ cover the quotient $\gag\backslash \T_{g, n}(\varepsilon)$. 
This provides that  $\gag\backslash \T_{g, n}(\varepsilon)$ is compact.

To prove that $\T_{g,n}(\varepsilon)$ is a real analytic manifold with corners,
we note by the collar lemma from hyperbolic geometry that when $\varepsilon$ is sufficiently small, for simple closed curves $\sigma_1,\dots, \sigma_m$ the level hypersurfaces 
$\{\ell_S(\sigma_1)=\varepsilon\},\dots,\{\ell_S(\sigma_m)=\varepsilon\}$
either have empty intersection or intersect transversally.
Since $\gag$ acts real analytically and properly on $\T_{g,n}(\varepsilon)$,
the quotient is a real analytic orbifold with corners, the desired conclusion.
\vspace{.1in}

Another approach for understanding the structure near infinity is to compactify the quotient $\gag\backslash \T_{g,n}$.  The approach involves including parameters for degenerate surfaces which are obtained by pinching geodesics. Let $\hat{\T}_{g,n}$ be the augmented Teichm\"{u}ller space, obtained by adding points describing marked stable Riemann surfaces of the same Euler number (see [Ab] for example).  The bordification $\hat\T_{g,n}$ is introduced by formally extending the range of the Fenchel-Nielsen coordinates.  For a length $\ell_S(\sigma)$ equal to zero, the twist is not defined and in place of the geodesic for $\sigma$ there appears a pair of cusps.  Following Abikoff, [Ab], and Bers [Be],  
the extended Fenchel-Nielsen coordinates describe marked possibly noded Riemann surfaces.  An equivalence relation is defined for marked noded Riemann surfaces and a construction is provided for adjoining to $\T_{g,n}$ frontier spaces (where subsets of lengths vanish) to obtain the augmented Teichm\"{u}ller space.  The space $\hat\T_{g,n}$ is not locally compact since in a neighborhood of $\ell_S(\sigma)$ vanishing, the angle $\theta_S(\sigma)$ has values filling $\mathbb R$. (Relatedly the action of $\gag$ on $\hat{\T}_{g,n}$ is not proper.)  Harvey considered a description of $\hat\T_{g,n}$ in terms of $SL(2,\mathbb R)$ representations and the Chabauty topology to show that $\gag\backslash \hat\T_{g,n}$ is compact [Hv2, Theorem 3.6.1].  As above, Proposition \ref{bers-region} also provides that the quotient is compact.   The Deligne-Mumford compactification provides that 
$\gag\backslash\T_{g,n}$ is contained in a compact orbifold \cite{dm}.  Bers considered families of Kleinian groups to also show that the quotient is contained in a compact orbifold, see [Be, \S 7].  
 
Now we provide another approach for Proposition \ref{finite-many}.

\begin{prop}\label{finite-many-1}
$\gag$ has finitely many conjugacy classes of finite subgroups.
\end{prop}
\begin{proof} As noted, $\gag$ acts properly on $\T_{g,n}$ and for a finite subgroup $H$ the set $\T_{g,n}^H$ of fixed points is nonempty.  The association of $\T_{g,n}^H$ to $H$ has the following properties.  The set $\T_{g,n}^H$ is a proper subset except for cases of generic involutions for the special topological types $(1,1),\,(0,4)$ and $(2,0)$.  For the generic point of $\T_{g,n}^H$ the full automorphism group of the Riemann surface is a realization of the group $H$ (modulo the ambiguity of extension by a special generic involution).  It follows that $\T_{g,n}^{H_1}$ is a $\gag$ translate of $\T_{g,n}^{H_2}$ if and only if $H_1$ and $H_2$ are conjugate subgroups (modulo the ambiguity of extension by the special generic involutions).  
As noted, the quotient $\gag\backslash\T_{g,n}$ is contained in a compact orbifold.  The orbifold locus (the image of the branching loci for the local manifold covers) locally has finitely many components.  By compactness the total number of orbifold components is finite.  In particular there are only finitely many orbifold components in $\gag\backslash \T_{g,n}$ and thus only finitely many distinct sets $\T_{g,n}^H$ modulo the action of $\gag$.   

\end{proof}

\section{Equivariant deformation retractions of Teichm\"{u}ller spaces}
 
We present the proof of Theorem \ref{main}.  Considerations begin with the systole function $\Lambda$ on $\T_{g,n}$, the length of the shortest closed geodesic.  We introduce an approximate gradient $V$ to the systole.  The vector field $V$ should define an equivariant deformation retraction of $\T_{g,n}$ onto the truncated Teichm\"{u}ller space $\T_{g,n}(\epsilon)$.  To this purpose $V$ is required to be continuous, $\gag$-invariant and to have support contained in the complement of a set with compact quotient.  The systole is the minimum of lengths $\Lambda=\min_{\sigma}\ell_S(\sigma)$.  At a point where multiple lengths $\ell_S(\hat\sigma)$ have the value $\Lambda$, the $V$-derivatives $V\ell_S(\hat\sigma)$ must agree for each length $\ell_S(\hat\sigma)=\Lambda$.  This compatibility condition is necessary for $\Lambda$ to have a continuous  $V$-derivative.  To understand the flow we add the inward pointing condition $V\Lambda\ge 0$ and require $V\Lambda=1$ on the complement of a set with compact quotient. 

We begin with gradients of small lengths.  The gradients of disjoint simple geodesics are always linearly independent [Wo4].  Precise Weil-Petersson gradient information is provided by the expansion
\begin{equation}
\label{WPprod}
\langle\grad\ell_S(\sigma)^{1/2},\grad\ell_S(\sigma')^{1/2}\rangle\ =\ \frac{\delta_{\sigma\sigma'}}{2\pi}\ +\ O((\ell_S(\sigma)\ell_S(\sigma'))^{3/2})
\end{equation}  
for $\sigma,\sigma'$ simple and either coinciding or disjoint; $\delta_{\sigma\sigma'}$ the Kronecker delta and the $O$-term constant depending only on an overall bound for the lengths 
[Wo5, Lemma 3.12].  Disjoint simple root-length gradients $(2\pi)^{1/2}\grad \ell_S(\sigma)^{1/2}$ 
are almost orthonormal.   The approach is to define for a set of geodesics $\mathcal S$ a vector field
\[
V\ =\ \sum\limits_{\sigma\in\mathcal S} \kappa_{\sigma}\,\grad \ell_S(\sigma)
\]
with functions $\kappa_{\sigma}$ determined by the system of equations
\begin{equation}
\label{kappaeqn}
V\ell_S(\sigma')\ = 1\quad \quad \mbox{for } \sigma'\in\mathcal S.
\end{equation}
The lengths $\ell_S(\sigma')$, $\sigma'\in\mathcal S$ increase at unit speed with respect to the $V$ flow.    

In [Har, \S 3, p. 172-173] the set $\mathcal S$ was taken as the totality of geodesics $\{\ell_S(\sigma)\le 3\epsilon \}$ for a suitable $\epsilon$.  It was claimed that each resulting function $\kappa_{\sigma}$ vanishes for $\ell_S(\sigma) > 2\epsilon$.  This is not true and consequently the sum over lengths $\{\ell_S(\sigma)\le 3\epsilon\}$ defines a vector field that is not continuous at a value $\ell_S(\sigma_0)=3\epsilon$ (for multiple small lengths a locus $\{\ell_S(\sigma_0)=3\epsilon\}$ intersects any sublevel set $\{\Lambda < \delta\}$).  The resulting flow is not continuous.

We now modify the approach to first locally define vector fields $V$.  The essential step is the definition of  open sets $U$ with sets of geodesics $\mathcal S$ providing nonsingular systems (\ref{kappaeqn}).  The open sets will satisfy a $\gag$ translation property and a condition on small lengths.  We work with the sublevel set $\{\Lambda <3\epsilon\}\subset \T_{g,n}$.  Each point of the sublevel set has a system of neighborhoods satisfying two conditions.  
The first condition is that an element of $\gag$ either stabilizes a neighborhood or translates the neighborhood to a disjoint set.  The existence of such systems of neighborhoods is a property of a discrete group action.  The second condition is that if $\ell_S(\sigma)=\Lambda$ somewhere on a neighborhood, then $\ell_S(\sigma)$ is bounded by $4\epsilon$ on the neighborhood.  The existence of such systems of neighborhoods follows from compactness of $\gag\backslash\hat\T_{g,n}$ or from general bounds for gradients of lengths [Wo5].  We write $\mathcal U$ for a neighborhood satisfying the two conditions.  
Define $\mathcal S_{\mathcal U}$ to be the set of geodesics $\sigma$ with $\ell_S(\sigma)=\Lambda$ somewhere on $\mathcal U$.   For $\epsilon$ sufficiently small, geodesics of length at most $4\epsilon$ are disjoint and the leading term of expansion (\ref{WPprod}) dominates. For $\epsilon$ sufficiently small, the equation (\ref{kappaeqn}) is nonsingular for each neighborhood.  We now write $V_{\mathcal U}$ for the solution of (\ref{kappaeqn}) for a neighborhood $\mathcal U$.   We note that the sets $\mathcal S_{\mathcal U}$ and vector fields $V_{\mathcal U}$ are canonically (without  choices) determined by the set $\mathcal U$.  If an element of $\gag$ stabilizes $\mathcal U$ then $V_{\mathcal U}$ is invariant with respect to the element.

We are ready to define a suitable open cover and a partition of unity.  For each neighborhood $\mathcal U$ introduce a smooth function $\psi_{\mathcal U}$, with $\supp(\psi_{\mathcal U})$ relatively compact in $\mathcal U$ and the supports also forming neighborhood systems at points.  The functions $\psi_{\mathcal U}$ are chosen to be invariant under the stabilizers $\stab(\mathcal U)\subset\gag$.  Now select a set of pairs $\{(\mathcal U_{\alpha},\psi_{\mathcal U_{\alpha}})\}_{\alpha\in\mathcal A}$  such that the supports $\{\supp(\psi_{\alpha})\}_{\alpha\in\mathcal A}$ provide a locally finite cover of  the locally compact set $\gag\backslash\{\Lambda\le 3\epsilon\}$.   Since as above an element of $\gag$ either stabilizes a pair or translates a pair to a second pair with  disjoint support, the set of pairs $\{(\supp(\psi_{\alpha})\circ\gamma^{-1},\psi_{\alpha}\circ\gamma)\}_{\alpha\in\mathcal A,\gamma\in\gag}$  provides a locally finite cover of $\{\Lambda\le 3\epsilon\}$ and is {\em invariant} in that an element of $\gag$ either stabilizes a pair or translates a pair to second pair in the set with disjoint support.   An overall consequence is that the function $\Psi=\sum_{\alpha\in\mathcal A,\gamma\in\gag}\psi_{\alpha}\circ\gamma\,$ is $\gag$-invariant and positive.  In particular the functions $\{\psi_{\alpha}\circ\gamma\slash\Psi\}_{\alpha\in\mathcal A,\gamma\in\gag}$ provide a $\gag$-invariant partition of unity of $\{\Lambda\le 3\epsilon\}$ (the partition support contains the sublevel set).  The set of triples $\{(\gamma^{-1}(\mathcal U_{\alpha}),V_{\alpha},\psi_{\alpha}\circ\gamma\slash\Psi)\}_{\alpha\in\mathcal A,\gamma\in\gag}$ is also invariant; an element of $\gag$ either stabilizes a triple or translates a triple to a second triple in the set with disjoint support.

Choose a non negative function $\phi(\ell)$ that is unity for $\ell\le 2\epsilon$ and vanishes for $\ell \ge 3\epsilon$.  The smooth vector field
\[
V\ =\ \phi(\Lambda)\sum\limits_{\alpha\in\mathcal A,\gamma\in\gag} (\gamma_*)(\psi_{\alpha}V_{\mathcal U_{\alpha}})
\]
satisfies $V\Lambda=1$ on $\{\Lambda\le 2\epsilon\}$, vanishes on $\{\Lambda\ge 3\epsilon\}$ and is $\gag$-invariant.   The time $\epsilon$ flow defines a smooth equivariant deformation retraction of $\T_{g,n}$ to the truncated Teichm\"{u}ller space $\T_{g,n}(\epsilon)$.

\begin{rem}
{\em An alternate approach for a deformation is as follows.
The augmented Teichm\"{u}ller space is the Weil-Petersson completion and a CAT(0) metric space [Wo1]  [DaW] [Ya].   The components of the bordification of $\hat\T_{g,n}$ are totally geodesic embeddings of products of lower dimensional Teichm\"{u}ller spaces. For a component boundary space $\T'$, determined by the free homotopy class $\sigma'$ represented by a node, there is a projection of $\hat\T_{g,n}$ to $\T'$ (the fibers of the projection are the geodesics realizing the distance to $\T'$).  The geodesics define a fibration of $\T_{g,n}$ with base space $\hat\T'$. A natural deformation from $\T_{g,n}$ to the submanifold  $\{x\in \T_{g,n}\mid \ell \ge \varepsilon\}$ is along fibers.  For multiple small lengths, fibrations might be combined to define a deformation with small lengths  increasing in a controlled manner.  
In [Wo5,  \S 4.2] it is shown that the fibers of the projections to boundary spaces are approximated to high order by integral curves of constant sums of 
gradients $\grad \ell(\sigma)^{1/2}$ for small lengths.  The descriptions by constant sums of root-length gradients and by fibers of projections essentially describe the same structure. From (\ref{kappaeqn}) the root-length gradient flow is essentially a reparameterization of the flow of $V$.  
}
\end{rem}

\noindent{\em Proof of Theorem \ref{contractible}.}  
Since $\T_{g,n}$ is contractible, Theorem \ref{main} provides that 
$\T_{g, n}(\varepsilon)$ is also contractible.
By the proof of Proposition \ref{finite-dim},
$\T_{g,n}^{H}$ is contractible; the equivariance of the deformation retraction
in Theorem \ref{main} provides that $(\T_{g, n}(\varepsilon))^{H}$
is a deformation retract of  $\T_{g,n}^{H}$ and hence contractible.

By Proposition \ref{finite-dim} again, $\T_{g,n}^{H}$ is non empty.
The above equivariant deformation retraction provides that $(\T_{g, n}(\varepsilon))^{H}$ is non empty.  Combined with Proposition \ref{cw-complex}, this provides that
$\T_{g,n}(\varepsilon)$ is a cofinite universal space for proper actions
of $\gag$.


\begin{thebibliography}{WX2}

\bibitem[Ab]{ab}W.Abikoff,  {\em The real analytic theory of Teichm\"{u}Ÿller space}, 
Lecture Notes in Mathematics, 820. Springer, 1980. vii+144 pp.

\bibitem[Ba]{ba}W.Baily, {\em 
On the moduli of Jacobian varieties}, Ann. of Math. (2) 71 (1960), 303--314. 

\bibitem[BCH]{bch}P.Baum, A.Connes,  N.Higson, {\em Classifying space for proper actions 
and $K$-theory of group $C^*$-algebras}, in {\em  $C^*$-algebras: 1943--1993},
pp. 240--291, Contemp. Math., 167, Amer. Math. Soc., Providence, RI, 1994.

\bibitem[Be]{be}L.Bers, {\em Spaces of degenerating {R}iemann surfaces},
in  {\em Discontinuous groups and Riemann surfaces (Proc. Conf., Univ.
             Maryland, College Park, Md., 1973)}, pp. 43--55. Ann. of Math. Studies, No. 79,
Princeton Univ. Press, 1974.

\bibitem[BS]{bs}A.Borel, J.P.Serre,  {\em Corners and arithmetic groups}, 
Comment. Math. Helv. 48 (1973) 436--491. 

\bibitem[Bri]{bri}M.Bridson,  {\em  Finiteness properties for subgroups of ${\rm GL}(n,\Z)$},
 Math. Ann. 317 (2000) 629--633.
 
\bibitem[Brou1]{brou1}S.Broughton,  {\em The equisymmetric stratification of the moduli space 
 and the Krull dimension of mapping class groups}, Topology Appl. 37 (1990) 101--113.
 
\bibitem[Brou2]{brou2}S.Broughton,  {\em  Normalizers and centralizers of elementary abelian subgroups of the mapping class group}, in {\em Topology '90 (Columbus, OH, 1990)}, 
 pp. 77--89, de Gruyter, Berlin, 1992.
 

\bibitem[Brow1]{brow1}K.Brown, {\em Finiteness properties of groups},
J. Pure Appl. Algebra 44 (1987), 45--75. 

\bibitem[Brow2]{brow2}K.Brown, {\em Cohomology of groups},
Graduate Texts in Mathematics, 87. Springer-Verlag,  1982. x+306 pp.

\bibitem[BHM]{bhm}M.B\"{o}kstedt, W.Hsiang, I.Madsen,  {\em
 The cyclotomic trace and algebraic $K$-theory of spaces}, Invent. Math. 111 (1993), 465--539.

\bibitem[Bu]{bu}P.Buser, {\em Geometry and spectra of compact Riemann surfaces},
Birkh\"{a}user, 1992.

\bibitem[DaW]{daw}G.Daskalopoulos, R.Wentworth,  {\em   Classification of {W}eil-{P}etersson isometries},
 125 (2003), {941--975},
 
\bibitem[DM]{dm}P.Deligne, D.Mumford, {\em The irreducibility of the space of curves of given genus}, Inst. Hautes \'Etudes Sci. Publ. Math.  36 (1969), 75--109.

\bibitem[FaM]{fam}B.Farb, D.Margalit, {\em A primer on mapping class groups}, preprint, 2005, 2008. 

\bibitem[Ga]{ga}F.Gardiner, {\em Teichm\"{u}ller theory and quadratic differentials},
 Wiley, New York, 1987. 

\bibitem[Gr]{gr}D.Grayson, 
{\em Reduction theory using semistability. II}, 
Comment. Math. Helv. 61 (1986), no. 4,
661--676.

\bibitem[Ham]{ham}U.Hamenst\"{a}dt, {\em Geometry of the mapping class groups I: 
Boundary amenability}, arXiv:math.GR/0510116, revised version, March 2008.

\bibitem[Har]{har}J.Harer,  {\em The cohomology of the moduli space of curves},
in {\em Theory of moduli},  pp. 138-221, Lecture Notes in Math., 1337, Springer,  1988. 

\bibitem[Hv1]{hv1}W.Harvey, {\em  Boundary structure of the modular group},
in {\em  Riemann surfaces and related topics},  pp. 245--251, Ann. of Math. Stud., 97, Princeton Univ. Press,  1981.

\bibitem[Hv2]{hv2}W.Harvey, {\em Spaces of discrete groups},
in {\em Discrete groups and automorphic functions ({P}roc. {C}onf.,
             {C}ambridge, 1975)}, pp.  {295--348}, {Academic Press}, {1977}. 

\bibitem[Il]{il}S.Illman, {\em Existence and uniqueness of equivariant triangulations of smooth 
proper $G$-manifolds with some applications to equivariant Whitehead torsion}, 
 J. Reine Angew. Math. 524 (2000), 129--183.
  

\bibitem[Iv1]{iv1}N.Ivanov, {\em  Mapping class groups}, in {\em Handbook of geometric topology}, 
pp. 523--633, North-Holland, Amsterdam, 2002.

\bibitem[Iv2]{iv2}N.Ivanov, {\em  Complexes of curves and Teichm\"{u}ller spaces},
Math. Notes 49 (1991), no. 5-6, 479--484

\bibitem[Iv3]{iv3}N.Ivanov, {\em 
Attaching corners to Teichm\"{u}ller space},
 Leningrad\ Math. J. 1 (1990), no. 5, 1177--1205.
 
\bibitem[Ji]{ji}L.Ji, {\em Integral Novikov conjectures and arithmetic groups containing torsion elements},  Comm. Anal. Geom. 15 (2007), no. 3, 509--533.
 
\bibitem[Ke]{ke}S.Kerckhoff, {\em The Nielsen realization problem}, 
Ann. of Math.  117 (1983), 235--265. 

\bibitem[Ki]{ki}Y.Kida, {\em The mapping class group from the viewpoint of measure equivalence theory}, arXiv:math/0512230. 

\bibitem[Le1]{le1}E.Leuzinger,
{\em An exhaustion of locally symmetric spaces by compact submanifolds 
with corners},   Invent. Math.  121  (1995),  no. 2, 389--410.

\bibitem[Le2]{le2}E.Leuzinger,
{\em On polyhedral retracts and compactifications of locally symmetric spaces},  Differential Geom. Appl.  20  (2004), no. 3, 293--318. 

\bibitem[Lo]{lo}S.Lojasiewicz, 
{\em Triangulation of semi-analytic sets}, Ann. Scuola Norm. Sup. Pisa (3) 18 (1964), 449--474.

\bibitem[Lu]{lu}W.L\"{u}ck,  {\em Survey on classifying spaces for families of subgroups},
in {\em Infinite groups: geometric, combinatorial and dynamical aspects}, pp.
269--322, Progr. Math., 248, Birkh\"{a}user, Basel, 2005. 

\bibitem[Mi]{mi}G.Mislin,  {\em Mapping class groups}, preprint, lecture notes of a talk
in M\"{u}nster, 2004.

\bibitem[Sa]{sa}L.Saper, {\em Tilings and finite energy retractions of locally symmetric spaces}, Comment. Math. Helv.  72  (1997),  no. 2, 167--202.

\bibitem[Se]{se}J.P.Serre, {\em Arithmetic groups}, in {\em  Homological group theory},   pp. 105--136, London Math. Soc. Lecture Note Ser., 36, Cambridge Univ. Press, 1979.

\bibitem[St]{st} P.Storm, {\em The {N}ovikov conjecture for mapping class groups
as a corollary of {H}amenst\"{a}dt's theorem}, arXiv:math/0504248v1.

\bibitem[Wo1]{wo1}S.Wolpert, {\em Geometry of the Weil-Petersson completion of Teichm\"{u}ller space}, in {\em Surveys in differential geometry}, Vol. VIII,   pp. 357--393, Surv. Differ. Geom., VIII, Int. Press,  2003.

\bibitem[Wo2]{wo2}S.Wolpert, {\em Geodesic length functions and the Nielsen problem},  
J. Differential Geom. 25 (1987), no. 2, 275--296. 

\bibitem[Wo3]{wo3}S.Wolpert, {\em Convexity of geodesic-length functions: a reprise}, 
in {\em  Spaces of Kleinian groups}, pp. 233--245, 
London Math. Soc. Lecture Note Ser., 329, Cambridge Univ. Press,  2006.

\bibitem[Wo4]{wo4}S.Wolpert, {\em The {F}enchel-{N}ielsen deformation},
 {Ann. of Math. (2)}, 115 (1982), {501--528}. 
   
\bibitem[Wo5]{wo5}S.Wolpert, {\em Behavior of geodesic-length functions on {T}eichm\"uller
             space}, {J. Differential Geom.},  79 (2008), {277--334}. 

\bibitem[Ya]{ya}{S.Yamada}, {\em On the geometry of {W}eil-{P}etersson completion of
             {T}eichm\"uller spaces}, Math. Res. Lett., 11 (2004), {327--344}.

\end{thebibliography}
\end{document}